\documentclass[12pt]{article}%
\usepackage{graphics}
\usepackage{color}
\usepackage{amsmath}
\usepackage{amsthm}
\usepackage{amssymb}
\usepackage{amsfonts}
\usepackage{graphicx}%
\newcommand{\be}{\mbox{{\bf E}}}
\newcommand{\bp}{\bar p}
\newcommand{\go}{\gamma_{0}}
\newcommand{\hn}{\hat N}
\newcommand{\hr}{\hat r}

\newcommand{\ou}{[0,1]}

\newcommand{\qn}{q_N\!}
\newcommand{\qnd}{q_N^2\!}
\newcommand{\sii}{\sigma_i}

\newcommand{\sln}{\sum_{l=1}^{n}}

\newcommand{\sk}{\mbox{SK}}

\newcommand{\ssn}{\Sigma_N}
\newcommand{\uu}{[-1,1]}
\newcommand{\1}{{\bf 1}}

\newcommand{\zn}{Z_N}

\newcommand{\cf}{{\cal F}}

\newcommand{\cl}{{\cal L}}

\newcommand{\cn}{{\cal N}}

\newcommand{\ga}{\gamma}
\newcommand{\gga}{\Gamma}

\newcommand{\ka}{\kappa}

\newcommand{\si}{\sigma}
\newcommand{\vp}{\varphi}

\newcommand{\E}{{\bf E}}

\newcommand{\R}{{\mathbb R}}
\newcommand{\Z}{\mathbb {Z}}

\newcommand{\lp}{\left(}
\newcommand{\rp}{\right)}
\newcommand{\lc}{\left[}
\newcommand{\rc}{\right]}

\newtheorem{theorem}{Theorem}[section]
\numberwithin{equation}{section}

\newtheorem{corollary}[theorem]{Corollary}

\newtheorem{hypothesis}[theorem]{Hypothesis}
\newtheorem{lemma}[theorem]{Lemma}

\newtheorem{proposition}[theorem]{Proposition}

\newtheorem{remark}[theorem]{Remark}

\begin{document}

\title{A central limit theorem for a localized version\\ of the SK model}
\author{S\'ergio de Carvalho  Bezerra \thanks{This author's research partially supported by CAPES.}
\hspace{0.3cm} Samy Tindel\\{\small \textit{Institut Elie Cartan,
Universit\'e de Nancy 1}}\\{\small \textit{BP 239,
54506-Vandoeuvre-l\`es-Nancy, France}}\\{\small
\texttt{[bezerra,tindel]@iecn.u-nancy.fr}}\vspace*{0.05in}}
\maketitle

\begin{abstract}
In this note, we consider a SK (Sherrington--Kirkpatrick)-type model on
$\mathbb{Z}^{d}$ for $d\geq1$, weighted by a function allowing to any single
spin to interact with a small proportion of the other ones.
In the thermodynamical limit, we investigate the equivalence of this
model with the usual SK spin system, through the study of the 
fluctuations of the free energy. 
\end{abstract}

\vspace{0.2in}\noindent\textbf{Key words and phrases:} spin glasses,
Sherrington-Kirkpatrick, localized mean-field model, cavity
method, stochastic calculus.

\vspace{0.3cm}


\section{Introduction}

This paper is concerned with a localized version of the 
Sherrington-Kirkpatrick model with external field, which
can be described in the following way:
for $N,d\geq1$, our space of configurations will be $\Sigma=\Sigma
_{N}=\{-1,1\}^{C_{N}}$, where $C_{N}$ is the finite lattice box $C_{N}%
=[-N;N]^{d}$ in $\mathbb{Z}^{d}$. For a given configuration $\sigma\in
\Sigma_{N}$, we will consider the Hamiltonian
\begin{equation}\label{introham}
-H_{N}\left(  \sigma\right)  =\frac{\beta}{\hat{N}^{d/2}}\sum_{(i,j)\in C_{N}%
}q\left(  \frac{i-j}{N}\right)  g_{(i,j)}\sigma_{i}\sigma_{j}+h\sum_{i\in
C_{N}}\sigma_{i}, 
\end{equation}
where $\beta$ stands for the inverse of the temperature of the system,
$\hat{N}=2N+1$, $(i,j)$ is the notation for a pair of sites $i,j\in
C_{N}$ (taken only once), $\left\{  g_{(i,j)}:(i,j)\in C_{N}\right\}  $ is a
family of {\small IID} standard centered Gaussian random variables, and $h$ represents
a constant positive external field, under which the spins tend to take the
value $+1$. Our localization is represented by the function $q$, which can be
thought of as a smooth frame, and which is only assumed to be defined on
$[-1,1]^{d}$ such that $q^{2}$ is of positive type, so that $q^{2}$ is
non-negative and invariant by symmetry about the origin. We also assume
$q^{2}$ is a continuous function, including at its periodic boundary
$-1\equiv 1$. 
The aim of our article is then to study the limit, when $N\to\infty$,
of the Gibbs measure $G_N(\si)$ defined on $\ssn$ by
$$
G_N(\si)=\frac{e^{-H_N(\si)}}{Z_N},
\quad\mbox{ where }\quad
Z_N=Z_N(\beta)=\sum_{\si\in\ssn}e^{-H_N(\si)},
$$
and more specifically, we will concentrate on 
the so-called free energy of the system, defined by:
\begin{equation}\label{def:free-energy}
p(\beta)=\lim_{N\to\infty}\be [p_N(\beta)]
=\mbox{a.\!\,s.}-\lim_{N\to\infty}p_N(\beta),
\quad\mbox{ where }\quad
p_N(\beta)=\frac{1}{\hn^d}\log(Z_N(\beta)).
\end{equation}

\vspace{0.3cm}

The model described by (\ref{introham}) can be considered as a finite range
approximation of the mean field SK model, 
associated with the Hamiltonian
$$
-\hat H_{N}\left(  \sigma\right)  
=\frac{\beta}{\hat{N}^{d/2}}\sum_{(i,j)\in C_{N}}
 g_{(i,j)}\sigma_{i}\sigma_{j}+h\sum_{i\in C_{N}}\sigma_{i},
$$
for which a large amount of information is now available \cite{GT1,Ta}.
It seems then natural to try to approximate the realistic spin glass system,
on which we have very little rigorous knowledge (see however \cite{NS}),
by our localized model (\ref{introham}), capturing some of the geometry
of the physical spin configuration, but still of a mean-field type
in the limit $N\to\infty$. One could then hope to perform an expansion
in $N$ in order to quantify the difference between the original
SK model and our model (\ref{introham}).

\vspace{0.3cm}

In fact, this kind of idea is not new, and goes back at least, in the
spin glass context, to \cite{FZ}. A version of our model with $h=0$
has been studied then in \cite{Toub}, and more recently, the Kac
limit of finite range spin glasses has been considered in 
\cite{GT,FT}. In these latter references, a slightly different point
of view is adopted: the finite range model depends on a given parameter
$\ga>0$, (which would be $1/N$ in our setting), and this localization 
parameter is sent to 0 \emph{after} the thermodynamical limit 
in $N$ is taken. It can be shown then, by some nice and soft interpolation
arguments, that in the limit $\ga\to 0$, the free energy of the localized
system is the same as the free energy of the SK model, for any value
of the parameter $\beta\ge 0$. Notice that the results contained
in \cite{GT,FT} cannot be applied directly to our model, since in our
case the limits $\ga\to 0$ and $N\to\infty$ are taken at the same time.
However, some slight modifications of the computations contained
in these papers would also show that our quantity $p_N(\beta)$
defined at (\ref{def:free-energy}) behaves like the free energy
of the SK model for large $N$.

\vspace{0.3cm}

The goal of our paper is then, in a sense, more modest than \cite{GT,FT},
since we will only deal with the high temperature region  of the model,
i.e. small values of $\beta$. On the other hand, our scope is to 
show that the equivalence between the SK model and our localized
model still holds, in the thermodynamical limit, for a second order
expansion of the free energy, that is in the central limit theorem
regime. More specifically, we will show the following limit result:
let $\go$ be the $L^2$-norm of $q$ in $\uu^d$. For $\beta,h>0$,
let also $s$ be the unique solution to the equation
\begin{align}\label{sdef}
s=\E \left[\tanh^2(\beta z\sqrt{s}+h)\right],
\end{align}
and $\sk(\beta,h)$ be the function
\begin{align}\label{Fdef}
\sk(\beta,h)=\beta^2(1-s)^2\slash
4+\log2+\E \left[\log\left[\cosh \left(\beta
z\sqrt{s}+h\right)\right]\right],
\end{align}
which represents the free energy of the SK model in the high temperature
region. Set also $p(\beta ,h)=\sk(\go^{1/2}\beta,h)$. Then, under
suitable conditions on $q$, we have
$$
(\cl)-\lim_{N\to\infty}\hn^{d/2}\lc p_N(\beta)- p(\beta ,h)\rc
=Y,
$$
where $Y$ is a centered Gaussian random variable with variance 
$\tau=\tau(\beta,h,\go)$. The announced equivalence, at the
CLT level, between the SK model and the localized one, springs then
from the fact that $\tau(\beta,h,1)$ is also
the variance of the Gaussian random variable which shows up in 
the central limit theorem of the SK case (see \cite{GT2,Ti}).

\vspace{0.3cm}

Let us say a few words about the method we have used in order to
get our result: since we are in the high temperature regime, we 
are allowed to use a cavity type method in order to compute the
limit of the overlap of the localized spin system. This yields
then the limit of the free energy in a straightforward manner.
It has also been shown in \cite{Ti} that the stochastic calculus
tools developed in \cite{CN} could be adapted  to the case
of spin glasses with external field. This induces a powerful method
for obtaining central limit theorems for the free energy, and 
interestingly enough, in this context, a dynamical point of
view gives some insight on a static stochastic problem.
We will
elaborate here on this method in order to treat the localized case,
by taking advantage systematically of the Fourier decomposition of $q$.

\vspace{0.3cm}

Our paper is divided as follows: at Section \ref{sec:simple-limit},
we compute the simple limit of the overlap function and of $p_N(\beta)$,
recovering the results obtained in \cite{FT} for the high temperature regime.
At Section \ref{sec:fluctuuations}, we derive the announced central limit 
theorem thanks to stochastic calculus tools.

\section{Simple limit of the free energy}\label{sec:simple-limit}

Recall that we are dealing with the system induced by the Hamiltonian
(\ref{introham}), and let us define some additional notations about
Gibbs averages: let $f:\ssn^n\to\R$ be a function of $n$ configurations,
with $n\ge 1$. Then we set
\begin{equation}\label{def:gibbs-av}
\rho(f)=\frac{1}{Z_N^n}\sum_{\si^1,\ldots,\si^n}
f(\si^1,\ldots,\si^n) \exp\lp -\sln H_N(\si^l) \rp,
\quad\mbox{ and }\quad
\nu(f)=\E\lc \rho(f) \rc.
\end{equation}
In the sequel of the paper, we will also  write $\partial_u\vp$
instead of $\partial\vp/\partial u$ for the derivative of
a function $\vp$ with respect to a parameter $u$, and
$\bp_N(\beta)=\E[p_N(\beta)]$.
With these notations in hand, our strategy in order to get the limit
of $p_N(\beta)$ will follow the classical steps of the cavity
procedure, namely:
\begin{enumerate}
\item
Find an expression for $\partial_\beta \bp_N(\beta)$ in terms of an
overlap-type function $R^{1,2}$.
\item
For a general function $f:\ssn^n\to\R$,
find a useful expression for $\partial_v\rho_v(f)$ along a suitable
path defined for $v\in\ou$, involving a Hamiltonian $H_{N,v}(\si)$.
\item
Compute $\rho_v(R^{1,2})$ inductively and deduce an expression for 
$\partial_\beta\bp_N(\beta)$, and then for $\bp_N(\beta)$.
\end{enumerate}
We will start by the first of these steps, for which we will introduce
a little more notation: first of all, we assume for the moment the
following basic hypothesis on $q$:
\begin{hypothesis}\label{hyp:q}
The function $q$ is continuous on $\uu^d$.
\end{hypothesis}
\noindent
In this case, 
the function $q^2:\uu^d\to\R_+$ can be
decomposed, as a function of $L^2(\uu^d)$, into a Fourier series
of the form:
\begin{equation}\label{dcp:q-fourier}
q^2(x)=\sum_{k\in\Z^d}\ga_k e^{\imath\pi k\cdot x},
\quad\mbox{ with }\quad
\sum_{k\in\Z^d}\ga_k^2<\infty
\quad\mbox{ and }\quad 
\gga\equiv\sum_{k\in\Z^d}\ga_k<\infty.
\end{equation}
Set also $R^{1,2}=\hn^{-d}\sum_{i\in C_N}\sii^1\sii^2$, and for any $k\in\Z^d$,
\begin{equation}\label{D: R_k}
R_k^{1,2}=\frac{1}{\hn^d}
\sum_{i\in C_N}e^{\frac{\imath\pi i\cdot k}{N}}\sii^1\sii^2,
\quad\mbox{ if $k\ne 0$,}\quad \mbox{ and }\quad
R_0^{1,2}= R^{1,2}-r,
\end{equation}
where $r$ is a positive constant, whose exact value will be determined
later on. Eventually, we will denote by $\qn$ the function 
$q(\cdot/N)$ defined on $[-N;N]^d$, and $\Z_*^d=\Z^d\backslash \{0\}$.
Then the following relation holds true:
\begin{proposition}\label{P:deriva}
For all $\beta>0$, we have
\begin{align*}
\partial_\beta\bp_{N}(\beta)  
=\frac{\beta}{\hat{N}^{2d}}\sum_{(i,j)\in C_{N}}
\qnd(i-j) 
-\frac{\beta}{2}\sum_{k\in\mathbb{Z}_*^{d}}\gamma_{k}
\mathbf{E}\lc \rho(R^{1,2}_k) \rc
-\frac{\beta}{2}\go \be\lc \rho(R^{1,2}) \rc
+\frac{\Gamma}{2\hat N^d}.
\end{align*}
\end{proposition}
\begin{proof}
Set $B(\sigma)=e^{-H_N(\sigma)}$. Then
the previous definitions and an elementary Gaussian integration
by parts yield (see \cite{Ta}):
\begin{align*}
\partial_\beta\bp_{N}(\beta)&=
\frac{\beta}{\hat N^{2d}}\mathbf{E}\left[\sum_{(i,j)\in C_N}q_N^2(i-j)-\sum_{(i,j)\in
    C_N}\sum_{\sigma^1,\sigma^2}
\frac{q_N^2(i-j)\sigma_i^1\sigma_j^1\sigma^2_i\sigma^2_j
B(\sigma^1)B(\sigma^2)}{Z_N^2} \right]\\
&=\frac{\beta}{\hat N^{2d}}
\lc
\sum_{(i,j)\in C_N}q_N^2(i-j)
-
\sum_{(i,j)\in C_N}q_N^2(i-j) 
\E\lc \rho(\sigma_i^1\sigma_j^1\sigma^2_i\sigma^2_j) \rc
\rc.
\end{align*}
Hence, using the decomposition (\ref{dcp:q-fourier}) of $q^2$, we obtain
\begin{align}
\partial_\beta\bp_{N}(\beta)&=\frac{\beta}{\hat
  N^{2d}}\mathbf{E}\left[\sum_{(i,j)\in C_N}q_N^2(i-j)-\sum_{(i,j)\in
    C_N}\rho\left(\sum_{k\in\mathbb{Z}^d}\gamma_ke^{\imath
      \pi\frac{i-j}{N}\cdot
      k}\sigma_i^1\sigma_j^1\sigma^2_i\sigma^2_j\right)
\right]\nonumber\\
&=\frac{\beta}{\hat
  N^{2d}}\mathbf{E}\Bigg[\sum_{(i,j)\in
    C_N}q_N^2(i-j)-\sum_{k\in\mathbb{Z}^d}\frac{\gamma_k}{2}\rho\left(\left|\sum_{i\in C_N}e^{\imath
      \pi\frac{i}{N}\cdot
      k}\sigma_i^1\sigma^2_i\right|^2\right)\nonumber\\&\hspace{6cm}+\sum_{k\in\mathbb{Z}^d}\frac{\gamma_k}{2}\rho\left(\sum_{i\in C_N}\left|e^{\imath
      \pi\frac{i}{N}\cdot
      k}\sigma_i^1\sigma^2_i\right|^2\right)
\Bigg],\nonumber
\end{align}
from which the desired expression is deduced easily.

\end{proof}

We are now ready to start the second step of the strategy mentioned above.

\subsection{The cavity method}

The cavity method will consist here in suppressing in a continuous way
the interactions between a certain site $m\in C_N$ and the remaining spins.
This will be done in the following way: set 
$\hat{C}_{N}^{m}=\{i\in C_{N};i\neq m\}$. We decompose, for all $\sigma
\in\Sigma_{N}$,%
\[
-H_{N}\left(  \sigma\right)  =-H_{N}^{m}\left(  \rho_{m}\right)  +\sigma
_{m}\left(  h+g_{m}\left(  \rho_{m}\right)  \right),
\]
where $\rho_{m}$ denotes the ordered spin values except for the $m$-th spin,
that is, the $(  \hat{N}^{d}-1) $-tuple 
$(  \sigma_{i};i\in\hat{C}_{N}^{m})$ and we denote
\[
g_{m}\left(  \rho_{m}\right)  =\frac{\beta}{\hat{N}^{d/2}}\sum_{i\in
\hat{C}_{N}^{m}}\qn(m-i)  g_{(m,i)}\sigma_{i},
\]
and we also use the notation
\[
-H_{N}^{m}\left(  \rho_{m}\right)  =\sum_{(i,j)\in\hat{C}_{N}^{m}}
\qn(i-j)  g_{\left(  i,j\right)  }\sigma_{i}\sigma_{j}+h\sum
_{i\in\hat{C}_{N}^{m}}\sigma_{i}.
\]
This new Hamiltonian is similar but not identical to the Hamiltonian
$-H_{N-1}$ on $\Sigma_{N-1}$.

\vspace{0.3cm}

The path we will build is now of the following form: consider a collection
$\{B_{i,m}(v);i\in \hat{C}_{N}^{m},v\in\ou \}$ of 
independent standard Brownian
motions, and $\{X(v);v\in\ou \}$ an independent  reversed time Brownian motion,
all defined on the probability space $(\Omega,\cf,P)$. Notice that
$g_{(i,m)}$ can be seen as the final value $B_{i,m}(1)$ 
of the Brownian motion $B_{i,m}$, and that $X$ is the solution
to a stochastic differential equation of the form
\begin{align}\label{eqrev}
X(v)=\eta-\int_0^v\frac{X(s)}{1-s}ds+W(v),\quad v\in[0,1],
\end{align}
where $\eta$ is a standard Gaussian random variable and $W$ another
Brownian motion, independent of the remainder of the randomness.
Notice that, for notational sake, we have written $B_{i,m}$
instead of $B_{(i,m)}$. Set then $\hr=\go^{1/2}r$ and
for $v\in\ou$, define
$$
-H_{N,v}\left(  \sigma\right)  =
-H_{N}^{m}\left(  \rho_{m}\right)  
+\frac{\beta}{\hn^{d/2}}\sum_{i\in\hat C_N^m} \qn(m-i) \sii \si_m B_{i,m}(v)
+\beta \hr^{1/2}\sum_{i\in\hat C_N^m} \sii X(t)
+ h \sum_{i\in\hat C_N^m} \sii.
$$
For $f:\ssn^n\to\R$, denote also by $\rho_v(f)$ the associated Gibbs average,
defined in a similar way to (\ref{def:gibbs-av}), and
$\nu_{m,v}(f)=\E[\rho_v(f)]$.

\vspace{0.3cm}

Recall now the following elementary lemma from \cite{Ti}:
\begin{lemma}\label{itorev}
For $k\ge 1$, let $\{B_l;\, l\le k  \}$ be a collection of independent
standard Brownian motions. Let also $X$ be the solution to
(\ref{eqrev}),  and $\vp:\R^{k+1}\to\R$ be a $C^2$ function having at
most exponential growth together with its first two derivatives.
Then, for any $v\in\ou$,
\begin{align*}
&\be\lc \vp(B_1(v),\ldots,B_k(v),X(v))  \rc\\
&=
\be\lc\vp(\eta)\rc 
+ \frac12\int_0^v\sum_{l=1}^{k}
\be\lc \partial_{x_lx_l}^{2}\vp(B_1(s),\ldots,B_k(s),X(s))  \rc ds\\
&\hspace{2cm}-\frac12\int_0^v 
\be\lc \partial_{x_{k+1}x_{k+1}}^{2}\vp(B_1(s),\ldots,B_k(s),X(s))  \rc ds.
\end{align*}
\end{lemma}

These preliminary tools yield the following differentiation rule:
\begin{proposition}\label{P:prin}
For any $f:\Sigma_{N}^{n}\rightarrow\R$ and 
$v\in\ou$, the derivative of $\nu_{m,v}(f)$ is given by
\begin{multline}\label{dnudt}
\partial_v\nu_{m,v}\left(  f\right)
=\beta^{2}\sum_{k\in\mathbb{Z}^{d}}
\gamma_{k}e^{\imath\pi m\cdot k/N}\Bigg(\sum_{1\leq
l<l^{\prime}\leq n}\nu_{m,v}\left(  f \bar R_{k}^{l,l^{\prime}}\sigma
_{m}^{l^{\prime}}\sigma_{m}^{l}\right) \\
-n\sum_{l=1}^{n}\nu_{m,v}\left(  f \bar R_{k}^{l,n+1}\sigma_{m}^{l}%
\sigma_{m}^{n+1}\right)  
+\frac{n(n+1)}{2}\nu_{m,v}\left(  f
\bar R_{k}^{n+1,n+2}\sigma_{m}^{n+1}\sigma_{m}^{n+2}\right)  \Bigg). 
\end{multline}
\end{proposition}

\begin{proof}
This result stems from an easy application of Lemma \ref{itorev}
to the function
$$
\vp\lp B_{i,m}(v),i\in C_N^m;\, X_m(v) \rp
\equiv \rho_{m,v}(f).
$$
The computations of the second derivatives of $\vp$ are a matter of easy
(though cumbersome) calculations, and are left to the reader for
sake of conciseness.

\end{proof}

Now that the variations of $\nu_{m,v}(f)$ have been computed, we can
proceed to get some bounds on the overlaps $R_k^{1,2}$.

\subsection{Bounds on the overlap}

Let us start with three lemmas
whose proofs follow essentially that of the corresponding
result in the cavity method for the standard Sherrington-Kirkpatrick model,
and which
will be combined with the explicit expression for $\nu_{m,v}$ in Proposition
\ref{P:prin} and a separate calculation of $\nu_{m,v}$ for $v=0$ to obtain
information on the actual expected overlaps, under $\nu=\nu_{m,1}$, i.e. for
$v=1$.

\begin{lemma}\label{nut<nu}
There exist two positive constants $c_{n,q,\beta}$ and
$c_{n,q,\beta}^{\prime}$ that depend only on $n,$ $q$ and $\beta$, and are
uniformly bounded in $\beta$ for $\beta\in\lbrack0,1]$, such that if $f$ is a
positive function on $\ssn^n$ then for all $m\in C_{N}$ ,
\begin{equation}
\nu_{m,v}\left(  f\right)  \leq c_{n,q,\beta}\,\nu\left(  f\right),
\label{nut<nu1}%
\end{equation}
and%
\begin{equation}
\left|  \nu_{m,v}\left(  f\right)  -\nu_{m,0}\left(  f\right)  \right|  \leq
c_{n,q,\beta}^{\prime}\beta^{2}\nu^{1/2}\left(  \left|  f\right|  ^{2}\right)
\left[  \sum_{k\in\mathbb{Z}^{d}}\gamma_{k}\nu^{1/2}\left(  \left|
R_{k}^{1,2}\right|  ^{2}\right)  \right]  . \label{nut<nu2}%
\end{equation}
\end{lemma}

\begin{proof}
See \cite[Propositions 2.4.6 and 2.4.7]{Ta}.

\end{proof}

For $v=0$, the expression of $\nu_v(f)$ can be simplified for a large class
of functions $f$ on $\ssn$. 

\begin{lemma}\label{nu0}
For fixed $m\in\Sigma_{N}$, let $f$ be a function on $\ssn^n$ 
that does not depend on the values $\sigma_{m}%
^{1},\sigma_{m}^{2},\ldots,\sigma_{m}^{n}$. Then for any subset $I$ of
$\left\{  1,\ldots,n\right\}  $ we have%
\[
\nu_{m,0}\left(  f\prod_{l\in I}\sigma_{m}^{l}\right)  =\mathbf{E}\left[
\tanh\left(  Y\right)  ^{|I|}\right]  \nu_{m,0}\left(  f\right),
\]
where $Y$ is the Gaussian random variable defined as:
\[
Y=\beta z\sqrt{ \gamma_{0}r}+h,
\]
with a standard normal variable $z$.
\end{lemma}

\begin{proof}
See \cite[Lemma 2.4.4]{Ta}.

\end{proof}

On the other hand, some symmetry properties for $v=1$ yield the
following kind of estimate:
\begin{lemma}\label{Overlapcor}
Set
\begin{equation}\label{def:delta-v}
\delta_v\equiv \sum_{k\in\Z^d}\ga_k \nu_v\lp |R_k^{1,2}|^2 \rp,
\quad\mbox{ for }\quad
v\in\ou.
\end{equation}
Then
$$
\delta_1=
\gamma_{0}\nu\left(  \left(  R^{1,2}-r\right)  \left(  \sigma
_{m}^{1}\sigma_{m}^{2}-r\right)  \right)  
+\frac{1}{\hat{N}^{d}}\sum_{k\in\mathbb{Z}_*^{d}}
\sum_{i\in C_{N}}\gamma_{k}
e^{\imath\pi i\cdot k/N}
\nu\left(\bar R_{k}^{1,2}\sigma_{i}^{1}\sigma_{i}^{2}\right).
$$
\end{lemma}

Now, in order to exploit Lemma \ref{nu0}, we must modify the above 
expression for
$\delta_1$ by completing the following two tasks:

\begin{itemize}
\item[{\bf (i)}] 
estimate the error made by replacing the arguments of $\nu_{i,0}$ by
functions that are of the same form as those in Lemma \ref{nu0};

\item[{\bf (ii)}] estimate the error made by replacing $\nu$ by $\nu_{m,0}$ (or
$\nu_{i,0}$ as appropriate).
\end{itemize}

\subsubsection{Task (i). Separation of cavity variable from others}

It is sufficient to replace $R^{1,2}$ by the same
quantity with the $m$-th term omitted: define
\[
 R_{-,m}^{1,2}:=\frac{1}{\hat{N}^{d}}\sum_{i\in\hat{C}_{N}%
^{m}}\sigma_{i}^{1}\sigma_{i}^{2}=R^{1,2}-\frac{1}{\hat{N}^{d}}\sigma_{m}%
^{1}\sigma_{m}^{2}=R^{1,2}+O\left(  \frac{1}{\hat{N}^{d}}\right),
\]
and similarly let%
\begin{align*}
R_{k,-,m}^{1,2}  &  :=\frac{1}{\hat{N}^{d}}\sum
_{i\in\hat{C}_{N}^{m}}\sigma_{i}^{1}\sigma_{i}^{2}e^{\imath\pi i\cdot k/N}\\
&  =R_{k}^{1,2}-\frac{1}{\hat{N}^{d}}\sigma_{m}^{1}\sigma_{m}^{2}e^{\imath\pi
m\cdot k/N}=R_{k}^{1,2}+O\left(  \frac{1}{\hat{N}^{d}}\right)  .
\end{align*}
Then we have the following relation, whose elementary proof
is omitted:
\begin{lemma}
\label{task1}%
For any $v\in\ou$ and $\delta_v$ defined at (\ref{def:delta-v}), it holds that:
\begin{align*}
\delta_{0}  &  =\gamma_{0}\nu_{m,0}\left(  
\left(  R_{-,m}^{1,2}-r\right)  
\left(  \sigma_{m}^{1}\sigma_{m}^{2}-r\right)
\right) \\
&  +\sum_{k\in\mathbb{Z}_*^{d}}\frac{\gamma_{k}}{\hat{N}^{d}}\sum_{i\in
C_{N}}\nu_{i,0}\left(  \bar R_{k,-,i}^{1,2}\sigma
_{i}^{1}\sigma_{i}^{2}\right)  e^{\imath\pi i\cdot k/N}+O\left(  \hat{N}%
^{-d}\right)  .
\end{align*}
\end{lemma}

\vspace{0.3cm}

Thanks to a Lemma by Lattala and Guerra, we can now choose
$r$ in order to eliminate one of the terms in our overlap calculation
$\delta_{0}$ (at $v=0$ with separated spins), as an immediate consequence
of Lemma \ref{nu0}.

\begin{lemma}
\label{lattalaguerra}For any choice of the parameters $\beta,\gamma_{0},h>0$,
the equation
\begin{equation}
r=\mathbf{E}\left[  \tanh^{2}\left(  \beta z\sqrt{\gamma_{0}r}+h\right)
\right]   \label{lattala}%
\end{equation}
has a unique solution $r\in\lbrack0,1]$, and we have $_{{}}$%
\[
\nu_{m,0}\left(  \left(   R_{-,m}^{1,2}-r\right)  \left(
\sigma_{m}^{1}\sigma_{m}^{2}-r\right)  \right)  =0.
\]
\end{lemma}

\vspace{0.3cm}

Let us take advantage of this relation, and try to write 
$\delta_{0}$ in terms of $r$: going back to the 
expression in Lemma \ref{task1}, we get
\begin{eqnarray}\label{ftvio}
\delta_0&=&
\sum_{k\in\Z_{*}^{d}}\frac{\gamma_{k}}{\hat{N}^{d}}\sum_{i\in C_{N}%
}e^{\imath\pi i\cdot k/N}\nu_{i,0}\left(  \bar R_{k,-,i}^{1,2}\sigma_{i}%
^{1}\sigma_{i}^{2}\right)  + O\lp  \frac{1}{N^d}\rp\nonumber\\
&=&r\sum_{k\in\Z_{*}^{d}}
\gamma_{k}\left(  A_{k}^{1}+A_{k}^{2}+A_{k}^{3}\right)
+ O\lp  \frac{1}{N^d}\rp,
\end{eqnarray}
with
\begin{align*}
A_{k}^{1} &  :=S_k \, \nu\left(  \bar R_{k}^{1,2}\right),
\qquad
A_{k}^{2}   :=\frac{1}{\hat{N}^{d}}\sum_{i\in C_{N}}\nu_{i,0}\left(
\bar R_{k,-,i}^{1,2}-\bar R_{k}^{1,2}\right)
e^{\imath\pi i\cdot k/N}\\
A_{k}^{3} &  :=\frac{1}{\hat{N}^{d}}\sum_{i\in C_{N}}\left[  \nu_{i,0}\left(
\bar R_{k}^{1,2}\right)  -\nu\left(  \bar R_{k}^{1,2}\right)
\right]  e^{\imath\pi i\cdot k/N},
\end{align*}
where we have set
\begin{equation}\label{def:ssk}
S_{k}=\frac{1}{\hat{N}^{d}}\sum_{i\in C_{N}}e^{\imath\pi i\cdot k/N}.
\end{equation}
We can bound now $\delta_0$ in the following way:
\begin{lemma}\label{bnd:delta0}
Recall that $\gga$ has been defined at (\ref{dcp:q-fourier}). Then
there exists a constant $\kappa$ that depends on $\beta$ and $q$ but is
bounded for $\beta$ bounded such that
\begin{equation}\label{A3}
|\delta_0|\le
\kappa \lp
\frac{\gga}{N}+
\beta^{2}\nu^{1/2}\left(  \left|  R_{k}^{1,2}\right|
^{2}\right)  \left[  \sum_{k\in\mathbb{Z}^{d}}\gamma_{k}\nu^{1/2}\left(
\left|  R_{k}^{1,2}\right|  ^{2}\right)  \right]  
\rp.
\end{equation}
\end{lemma}

\begin{proof} 
Go back to relation (\ref{ftvio}), and let us bound the terms
$A_1^k, A_2^k, A_3^k$ for $k\in\Z^d$. First of all, the estimation
of $A_1^k$ is controlled by $S_k$. However,
$$
S_k
=
\frac{1}{\hn^d}\sum_{j_1=-N}^{N}\cdots \sum_{j_d=-N}^{N}
\prod_{l=1}^{d}e^{\imath\pi k_l j_l/N}
=\frac{1}{\hn^d} \prod_{l=1}^{d} 
\lc (-1)^{k_l}\1_{(k_l\ne 0)}+\hn \1_{(k_l= 0)}\rc,
$$
by an elementary argument on sums of geometric sequences. Hence,
\begin{equation}\label{bnd:ssk}
|S_k|\le \frac{\kappa}{N^{Z(k)}},
\quad\mbox{ where }\quad
Z(k)=
\mbox{number of components of $k$ that are non-zero},
\end{equation}
and thus,
\begin{equation}\label{bnd:sum-ak1}
\Big|
\sum_{k\in\Z_{*}^{d}}\gamma_{k}  A_{k}^{1}
\Big|
=
\Big|
\nu\lp \bar R_k^{1,2} \rp
\sum_{k\in\Z_{*}^{d}}\gamma_{k}  S_k
\Big|
\le \frac{\kappa\gga}{N}.
\end{equation}
Furthermore, it is easily checked that
\begin{equation}\label{a2-a3}
|A_k^2|\le \frac{\ka}{N}
\quad\mbox{ and }\quad
|A_k^3|\le
\kappa 
\beta^{2}\nu^{1/2}\left(  \left|  R_{k}^{1,2}\right|
^{2}\right)  \left[  \sum_{k\in\mathbb{Z}^{d}}\gamma_{k}\nu^{1/2}\left(
\left|  R_{k}^{1,2}\right|  ^{2}\right)  \right],
\end{equation}
by applying inequality (\ref{nut<nu2}) in Lemma \ref{nut<nu}.
Our claim is then proved easily by putting together 
(\ref{bnd:sum-ak1}) and (\ref{a2-a3}).

\end{proof}

\begin{remark}
In the remainder of the article, $\kappa$ will stand for a positive constant
depending on $\beta$ and $q$ and that is uniformly bounded in the range of our
parameter $\beta$; we will allow $\kappa$ to change from line to line
\end{remark}

\vspace{0.3cm}

Notice that, for our result on the fluctuations of $Z_N$, we will need an 
improved
bound on $A_k^1$. This can be achieved under the following additional
condition:
\begin{hypothesis}\label{hyp:q2}
Going back to the decomposition (\ref{dcp:q-fourier}), we assume
that there exists an integer $\hat d > d\slash 2$ such that, for every
$k\in\mathbb{Z}_*^d$ such that $\gamma_k\neq 0$, the number $Z(k)$ of
components of $k$ that are non-zero satisfies $Z(k)\geq \hat d$.
\end{hypothesis}
Then Lemma \ref{bnd:delta0} can be enhanced in the following way:
\begin{corollary}\label{bnd:delta0-bis}
Assume $q$ satisfies Hypothesis \ref{hyp:q} and \ref{hyp:q2}. Then
$$
|\delta_0|\le
\kappa \lp
\frac{\gga}{N^{\hat d}}+
\beta^{2}\nu^{1/2}\left(  \left|  R_{k}^{1,2}\right|
^{2}\right)  \left[  \sum_{k\in\mathbb{Z}^{d}}\gamma_{k}\nu^{1/2}\left(
\left|  R_{k}^{1,2}\right|  ^{2}\right)  \right]  
\rp.
$$
\end{corollary}

\begin{proof}
This is trivially checked by going through the computations of
Lemma \ref{bnd:delta0} again, and taking into account (\ref{bnd:ssk}).

\end{proof}

\subsubsection{Task (ii). Difference between the overlaps at
$v=0$ and $v=1$}

We are ready to state and prove the result which completes task (ii).
\begin{lemma}
\label{task2}There exists a constant $\kappa$ that depends on $\beta$ and $q$
but is bounded for $\beta$ bounded, such that%
\begin{equation}
\left|  \delta_1-\delta_0\right|  \leq  
+\kappa\beta^{2}\gga  
\sum_{k\in\mathbb{Z}^{d}}\gamma_{k}\nu\left(  \left|  R_{k}^{1,2}\right|
^{2}\right) .  \label{O-O0}%
\end{equation}
\end{lemma}

\begin{proof} 
We can first write, using Lemma \ref{nut<nu},
\begin{align}\label{bnd:zd-star}
&  \sum_{k\in\Z_{*}^{d}}\gamma_{k}\left|  \frac{1}{\hat{N}^{d}}%
\sum_{i\in C_{N}}\left[  \nu\left(  \bar R_{k}^{1,2}\sigma_{i}^{1}%
\sigma_{i}^{2}\right)  -\nu_{i,0}\left(  \bar R_{k}^{1,2}\sigma_{i}%
^{1}\sigma_{i}^{2}\right)  \right]  e^{\imath\pi i\cdot k/N}\right| 
\nonumber\\
&  \leq\kappa\beta^{2}\sum_{k\in\mathbb{Z}_*^{d}}
\gamma_{k}\nu^{1/2}\left(  \left|  R_{k}^{1,2}\right|  ^{2}\right)  \left[
\sum_{k\in\mathbb{Z}^{d}}\gamma_{k}\nu^{1/2}\left(  \left|  R_{k}%
^{1,2}\right|  ^{2}\right)  \right]  .
\end{align}
Similarly we can obtain%
\begin{align}\label{bnd:zd}
&  \gamma_{0}\left|  \nu\left(  \left(  R^{1,2}-r\right)  \left(  \sigma
_{m}^{1}\sigma_{m}^{2}-r\right)  \right)  -\nu_{i,0}\left(  \left(
R^{1,2}-r\right)  \left(  \sigma_{m}^{1}\sigma_{m}^{2}-r\right)  \right)
\right| \nonumber\\
&  \leq\gamma_{0}\kappa\beta^{2}\nu^{1/2}\left(  \left|  R_{0}^{1,2}\right|
^{2}\right)  \left[  \sum_{k\in\mathbb{Z}^{d}}\gamma_{k}\nu^{1/2}\left(
\left|  R_{k}^{1,2}\right|  ^{2}\right)  \right]  .
\end{align}
We now get, by putting together (\ref{bnd:zd-star}) and
(\ref{bnd:zd}) and applying Jensen's inequality, that
\begin{align*}
\left|  \delta_1-\delta_0\right|   &  =   \kappa\beta^{2}\left[  \sum_{k\in
\mathbb{Z}^{d}}\gamma_{k}\nu^{1/2}\left(  \left|  R_{k}^{1,2}\right|
^{2}\right)  \right]  ^{2}
  \leq    \kappa\beta^{2}\gga  
\sum_{k\in\mathbb{Z}^{d}%
}\gamma_{k}\nu\left(  \left|  R_{k}^{1,2}\right|  ^{2}\right)  ,
\end{align*}
which proves the lemma.

\end{proof}

\subsubsection{Self-averaging overlap limit}

We are now in a position to estimate $\delta_1$. We show that for small
$\beta$, this expected total overlap, which is recentered using the value $r$,
converges to $0$ at the speed $1/N$ as long as $q$ is a continuous function.
\begin{proposition}
\label{prop1}Let $\beta>0$ and let $r=r\left(  \beta\right)  $ be the solution
of (\ref{lattala}). Let $\kappa$ be the constant defined in 
Lemma \ref{task1} and Lemma \ref{task2},
i.e. $\kappa$ is a constant that depends on $\beta$ and $q$ but is bounded for
$\beta$ bounded. Assume that $q$ satisfies Hypothesis \ref{hyp:q},
and that $\beta$ is so small that $2\kappa\beta^{2}\gga <1$.
In that case, we have, for $N$ large enough,
with $R_{k}^{1,2}$ defined by relation (\ref{D: R_k}),
\[
0\leq\nu\left(  \sum_{k\in\mathbb{Z}^{d}}\gamma_{k}\left|  R_{k}^{1,2}\right|
^{2}\right)  \leq
\frac{r\gga}{(1-2\ka\beta^2\gga)N}.
\]
\end{proposition}

\begin{proof} 
We have, using (\ref{A3}) and (\ref{O-O0}),
\begin{align}\label{finalestimateforcurlyO}%
0 &  \leq\delta_1=\nu\left(  \sum_{k\in\mathbb{Z}^{d}}\gamma_{k}\left|
R_{k}^{1,2}\right|  ^{2}\right)  \leq
\left|  \delta_0-\delta_1\right|
+r |\delta_0|\nonumber\\&\le 
2\kappa\beta^{2}\gga 
\nu\left(  \sum_{k\in\mathbb{Z}^{d}}\gamma_{k}\left|  R_{k}^{1,2}\right|
^{2}\right)+ O\lp \frac{1}{N} \rp,
\end{align}
where we recall that $O(N^{-1})$ is a function that tends
 to zero as fast as $N^{-1}$ and that 
this convergence holds uniformly in all parameters. Moreover,
since $\kappa$ is bounded for
$\beta$ bounded, for $\beta$ sufficiently small we can make 
$2\kappa\beta^{2}\gga$ smaller than $1$.
The result of the proposition follows.

\end{proof}

\begin{corollary}
\label{cor1}Under the same assumptions as in Proposition \ref{prop1}, 
but assuming additionally that condition \ref{hyp:q2}
holds true, then the conclusion of Proposition
\ref{prop1} holds with $N$ replaced by $N^{\hat d }$.
\end{corollary}

\begin{proof} 
This follows trivially from the proof of Proposition
\ref{prop1} if we modify the argument in order to take into account
Corollary \ref{bnd:delta0-bis}.

\end{proof}

\subsection{Consequence for the partition function}\label{limmean}

We can now apply the previous computations in order to get the simple
limit of $\bp_N(\beta)$, which recovers, in the high temperature
region, the results contained in \cite{FT}.
\begin{theorem}
\label{TheopN}Under the hypotheses of Proposition \ref{prop1}, we have%
\begin{equation}
\left| \bp_{N}\left(  \beta\right)  -\sk\left(  \gamma_{0}^{1/2}
\beta,h\right)  \right|  \leq\frac{C\left(  \beta\right)
}{N}, \label{limitofpN}%
\end{equation}
where the constant $C$ depends on $h$, $q$, and $\beta$, and is bounded for
$\beta\in\lbrack0,\beta_{0}]$.
\end{theorem}

\begin{proof} 
Recall from Proposition \ref{P:deriva} that
\begin{equation}
\partial_\beta \bp_{N}=\frac{\beta}{\hat{N}^{2d}}\sum_{(i,j)\in
C_{N}}\qnd(i-j) -\frac{\beta}{2}\sum_{k\in\Z_*^d}
\gamma_{k}\nu\left(  \left|  R_{k}^{1,2}\right|  ^{2}\right)
-\frac{\beta}{2}\gamma_{0}\nu\left(  \left|  R_{0}^{1,2}+r\right|
  ^{2}\right)+\frac{\Gamma}{2\hat N^d}.
\label{whereB1}%
\end{equation}
The first term on the right-hand side can be handled easily: indeed, we
have
\begin{eqnarray*}
\frac{1}{\hat{N}^{2d}}\sum_{(i,j)\in C_{N}}\qnd(i-j)
&=&
\frac{1}{\hat{N}^{2d}}\sum_{(i,j)\in C_{N}}\sum_{k\in\Z^d}\gamma_k
e^{\imath\pi k\cdot i/N} e^{-\imath\pi k\cdot j/N}\\
&=&
\frac{1}{2\hat{N}^{2d}}\sum_{i\ne j}\sum_{k\in\Z^d}
\gamma_ke^{\imath\pi k\cdot i/N} e^{-\imath\pi k\cdot j/N}\\
&=&\frac12\lp \go+ \sum_{k\in\Z_*^d}\ga_k| S_k|^2 -\frac{\gga}{\hn^d}\rp,
\end{eqnarray*}
where $S_k$ has been defined at (\ref{def:ssk}). Thus,
$$
\Big|
\frac{\beta}{\hat{N}^{2d}}\sum_{(i,j)\in C_{N}}\qnd(i-j)
-
\frac{\beta\go}{2}
\Big|
\le \frac{\ka}{N},
$$
and Proposition
\ref{prop1} then easily yields, for $\beta<\beta_{0}$,%
\begin{equation}\label{delpestim}
\left|  \partial_\beta \bp_{N}-\frac{\beta\go}{2}-\frac{r^{2}\gamma
_{0}\beta}{2}\right|  \leq
\frac{\ka}{N}
+\gamma_{0}r\beta\left|  \nu\left(  R_{0}^{1,2}\right)  \right|.
\end{equation}
Furthermore, it can be shown, along the same lines as in \cite{Ta},
that for $\beta<\beta_0$, we have 
$$
\left|  \nu\left(  R_{0}^{1,2}\right)
\right|  \leq K\left(  \beta,q\right)  /N,
$$ 
where $K\left(  \beta,q\right)$ is bounded for $\beta$ bounded.
This inequality and (\ref{delpestim}) now imply
\[
\left|  \partial_\beta \bp_{N}-\frac{\gamma_{0}\beta}{2}\left(
1-r^{2}\right)  \right|  \leq\frac{\kappa}{N}%
\]
where $\kappa$ depends only on $\beta,h,q$ and is bounded for $\beta\in
\lbrack0,\beta_{0}]$. The theorem follows by integrating 
$\partial_\beta \bp_{N}$, and using trivial calculations 
and known facts about the
function $\sk$.
\end{proof}

\vspace{0.3cm}

The final result we present in this section shows that while the complete
structure of $q$ does not seem to effect the limiting behavior of the
partition function beyond the average value $\gamma_{0}$ of $q^{2}$, the speed
of convergence towards this value may depend heavily on the behavior of $q$.
We show that the speed can be increased to the order $N^{-\hat d }$ as long
as the hypotheses of Corollary \ref{cor1} hold. 

\begin{corollary}
Under the hypotheses of Corollary \ref{cor1}, we have%
\[
\left| \bp_{N}\left(  \beta\right)  -\sk\left(  \go^{1/2}\beta,h\right)  \right|  \leq\frac{\beta C\left(  \beta\right)
}{N^{\hat d }}%
\]
where the constant $C$ depends on $h$, $q$, and $\beta$, and is bounded for
$\beta\in\lbrack0,\beta_{0}]$.
\end{corollary}

\begin{proof} 
In the proof of Theorem \ref{TheopN}, 
some estimates are already of order $N^{-d}$. 
For the others, we may use Corollary \ref{cor1}
instead of Proposition \ref{prop1} in all its occurrences. This improves all
estimates that were originally of order $N^{-1}$ to the order $N^{-\hat d %
}$, with the exception of the estimation of the first term on the right-hand
side of (\ref{whereB1}). But here again, one can easily check that
this term is of order $N^{-\hat d}$ under the conditions of
Corollary \ref{cor1}.

\end{proof}

\section{Fluctuations of the free energy}\label{sec:fluctuuations}

In this section, we will turn to our main aim, that is the central
limit theorem governing the fluctuations of $Z_N$. More specifically, we 
will get the following:
\begin{theorem}\label{T:tcl}
For $\beta$ small enough, $t\in[0,1]$, if $q$ satisfies Hypothesis
\ref{hyp:q} and \ref{hyp:q2}, then
$$
(\cl)-\lim_{N\to\infty}\hn^{d/2}\lc p_N(\beta)-p(\beta,h)\rc
=Y,
$$
where $Y$ is a centered Gaussian random variable with variance 
$\tau$, and $\tau$ is given by
\begin{equation}\label{def:tau}
\tau=\hat\tau
-\frac{\beta^2\gamma_0r^2}{2},
\quad\mbox{ where }\quad
\hat\tau=\mathrm{Var}\lc \log(\cosh(\beta\sqrt{\gamma_0r}z+h))\rc
\end{equation}
with $r$ defined by (\ref{lattala}), and a standard Gaussian random variable
$z$.
\end{theorem}

This kind of result is usually obtained by letting all the interactions
between spins tend to 0 at once, and this procedure can be somewhat 
simplified by considering the computations from a stochastic
calculus point of view. This leads us to consider a new path
$t\mapsto H_{N,t}(\si)$ defined for $t\in\ou$ by
\begin{align}\label{E:Ham1}
-H_{N,t}(\sigma)=\frac{\beta}{\hat{N}^{d\slash 2}}\sum_{(i,j)\in
C_N}B_{i,j}(t)\qn(i-j)\sigma_i\sigma_j+\beta
\hat r^{1\slash 2}\sum_{i\in C_N}X_i(t)\sigma_i+\sum_{i\in
C_N}h\sigma_i,
\end{align} 
where $\hr=\gamma_0r$ and $\{B_{i,j};(i,j)\in C_N\}$ is again a collection of independent
standard Brownian motions, and $\{X_i;i\in C_N  \}$ is a family of
independent reversed time Brownian motion, which can be seen as 
the solution to some stochastic differential equations of the
form (\ref{eqrev}), for a family $\{\eta_i;i\in C_N  \}$
(resp. $\{W_i;i\in C_N  \}$) of standard Gaussian random variables
(resp. of independent Brownian motions). Here again, we will assume
that all these objects are defined on the same probability space
$(\Omega,\cf,P)$. Notice that, in order to spare notations,
we have called this modified Hamiltonian $H_{N,t}$ again,
like in Section \ref{sec:simple-limit}, hoping that this won't
lead to any confusion. We will also denote by $\zn(t)$ and
$\rho_t(f)$ the partition function and the Gibbs average
associated with $H_{N,t}$. Eventually, for $t\in\ou$, we define
$$
p_{\beta,h,t}=\frac{\beta^2\gamma_0 t}{4}\left(1-r\right)^2
+\log(2)
+\mathbf{E}\lc\log(\cosh(\beta\sqrt{\gamma_0r}z+h))\rc,
$$
and we observe that $p_{\beta,h,1}=p(\beta,h)$.
Once this path has been defined, the strategy
which leads to Theorem \ref{T:tcl} can be summarized as follows: 
apply It\^o's formula in order to:
\begin{enumerate}
\item
Compute the variations of $t\mapsto e^{-H_{N,t}(\si)}$.
\item
Get an equation for the evolution of 
$Y_N(t)\equiv \hn^{d/2}[p_N(\beta)-p_{\beta,h,t}]$.
\item
Find a limit for $\be[e^{\imath u Y_N(t)}]$ for any $u\in\R$.
\end{enumerate}
We will detail now this global strategy.

\subsection{Preliminary computations}

Let us first study the dynamics of $t\mapsto e^{-H_{N,t}(\si)}$:
\begin{proposition}\label{P:exp}For $t\in[0,1]$ and $\sigma\in \Sigma_N$,
we have
$$
e^{-H_{N,t}(\sigma)}=\exp\left(\sum_{i\in C_N}\sigma_i(\beta \hat
r^{1\slash 2}\eta_i+h)\right) + F_N^1(t)+F_N^2(t),
$$
with
\begin{eqnarray*}
F_N^1(t)&=&
\frac{\beta}{\hat{N}^{d\slash 2}}\sum_{(i,j)\in C_N}
\qn(i-j)\sigma_i\sigma_j\int_0^te^{-H_{N,s}(\sigma)}dB_{i,j}(s)\\
&&+\beta \hat r^{1\slash 2}\sum_{i\in C_N}
\sigma_i\int_0^te^{-H_{N,s}(\sigma)}dW_i(s),
\end{eqnarray*}
and
\begin{eqnarray*}
F_N^2(t)&=&
-\beta \hat r^{1\slash 2}\sum_{i\in C_N}
\sigma_i\int_0^te^{-H_{N,s}(\sigma)}\frac{X_i}{1-s}ds\\
&&+\frac{\beta^2}{2}\left(\frac{1}{\hat{N}^d}\sum_{(i,j)\in C_N}\qnd(i-j)
+\hat r\hat N^d\right)\int_0^te^{-H_{N,s}(\sigma)}ds.
\end{eqnarray*}
\end{proposition}
\begin{proof}
The exponential function being a $C^2$ function with a nicely
controlled growth, we can apply Ito's formula and we obtain
\begin{align*}
e^{-H_{N,t}(\sigma)}
&=e^{-H_{N,0}(\sigma)}-\int_0^te^{-H_{N,s}(\sigma)} dH_{N,s}(\sigma)
+\frac{1}{2}\int_0^te^{-H_{N,s}(\sigma)}d\langle H_{N,\cdot}\rangle_s,
\end{align*}
where $\langle M\rangle_t$ denotes the quadratic variation process
of a semi-martingale $M$.
Now, we can evaluate the quantity $\langle H_{N,\cdot}\rangle _s$, since 
all of the
$B_{i,j}$, $W_i$ are independent and using the fact that any finite
variation process have a null quadratic variation. We get
\begin{eqnarray*}
\langle H_{N,\cdot}\rangle _t&=&\frac{\beta^2}{\hat N^d}
\sum_{(i,j)\in  C_N}
\qnd(i-j)(\sigma_i\sigma_j)^2t
+\beta^2\hat r\sum_{i\in C_N}\sigma_i^2t\\
&=& \frac{\beta^2}{\hat N^d}
\sum_{(i,j)\in  C_N}
\qnd(i-j) t
+\beta^2\hat r \hn^d t,
\end{eqnarray*}
from which our claim is easily shown.

\end{proof}

\vspace{0.3cm}

We will turn now to the fluctuations of $p_N(\beta)$. Namely set,
for $t\in\ou$,
\begin{align}
Y_N(t)=\hat N^{d\slash2}
\left(\frac{1}{\hn^d}\log(\zn(t))-p_{\beta,h,t}\right),
\end{align}
and define also the function $\Phi:\R\to\R$ by
\begin{align*}
\Phi(x)=\log(\cosh (\beta\hat r^{1\slash 2}x+h)).
\end{align*}
Then the semi-martingale $Y_N$ can be decomposed in the following way:
\begin{proposition} 
Recall that $\gga=\sum_{k\in\Z^d}\ga_k$. Then,
for $t\in\ou$, $Y_N(t)$ satisfies:
\begin{align*}
Y_N(t)=U_N+\sum_{l\leq 2}
M_{l,N}(t)-(V_{1,N}(t)-V_{2,N}(t))+V_{3,N}(t),
\end{align*}
where the random variable $U_N$ and the processes $M_{l,N}$ and
$V_{k,N}$ are defined by:
\begin{align*}
U_N&= \hat N^{d\slash 2}\left(\frac{1}{\hat N^d}\sum_{i\in C_N}\Phi(\eta_i)-\E \Phi(z)  \right) \\
M_{1,N}(t)&=\frac{\beta}{\hat N^{d}}\sum_{(i,j)\in
C_N}\qn(i-j)\int_0^t\rho_s(\sigma_i\sigma_j)dB_{i,j}(s)\\
M_{2,N}(t)&=\frac{\beta\hat r^{1\slash 2}}{\hat N^{d\slash 2}}
\sum_{i\in C_N}\int_0^t\rho_s(\sigma_i)dW_i(s) \\
V_{1,N}(t)&=\frac{\beta\hat r^{1\slash 2}}{\hat N^{d \slash
2}}\sum_{i\in
C_N}\int_0^t\rho_s(\sigma_i)\frac{X_i(s)}{1-s}ds   \\
V_{2,N}(t)&=\frac{\beta^2\hat r}{\hat N^{d\slash 2}}\sum_{i\in C_N}\int_0^t\rho_s(1-\sigma_i^1\sigma_i^2)ds   \\
V_{3,N}(t)&=
\lc  \frac{\gga}{\hn^{d/2}}
+\frac{2}{\hat N^{3d\slash 2}}\sum_{(i,j)\in C_N}\qnd(i-j)
- \gamma_0\hn^{d/2}\rc
\frac{\beta^2 t}{4}
\\
&-\frac{\hat N^{d\slash 2}\beta^2}{4}
\sum_{k\in\mathbb{Z}^d}\gamma_k \int_0^t\rho_s\left(\left|R_k^{1,2}\right|^2\right)ds.
\end{align*}
\end{proposition}

\begin{proof}
Notice that $\zn(t)$ is almost surely a strictly positive random variable.
Thus, It\^o's formula can be applied to $\log(\zn(t))$, and we obtain
\begin{align}\label{ito:log-zn}
\log(Z_N(t))=\log(Z_N(0))
+\int_0^t\frac{dZ_N(s)}{Z_N(s)}
-\frac{1}{2}\int_0^t\frac{d\langle Z_N\rangle _s}{Z_N^2(s)}.
\end{align}
The first two terms in the right hand side of (\ref{ito:log-zn})
are easily computed. Indeed, it is easily checked that
\begin{align}\label{E:zn0}
\log(Z_N(0))=\hn^d\log2 +\sum_{i\in C_N}\log[\cosh (\beta \hat
r^{1\slash 2}\eta_i+h)].
\end{align}
Furthermore, invoking the fact that 
$\zn(t)=\sum_{\si\in\ssn}e^{-H_{N,t}(\si)}$ and Proposition \ref{P:exp},
we have
\begin{multline}\label{E:rec1}
\int_0^t\frac{dZ_N(s)}{Z_N(s)}=\frac{\beta}{\hat N^{d\slash
2}}\sum_{(i,j)\in
C_N}\qn(i-j)\int_0^t\rho_s(\sigma_i\sigma_j)dB_{i,j}(s)
+\beta
\hat r^{1\slash 2}\sum_{i\in
C_N}\int_0^t\rho_s(\sigma_i)dW_i(s)\\
-\beta \hat r^{1\slash 2}\sum_{i\in
C_N}\int_0^t\rho_s(\sigma_i)\frac{X_i(s)}{1-s}ds
+\left(\frac{\beta^2}{2\hat N^d}\sum_{(i,j)\in
C_N}\qnd(i-j)+\frac{\beta^2\hat r\hat
N^d}{2}\right)t.
\end{multline}
Thus, putting together (\ref{ito:log-zn}), (\ref{E:zn0}) and (\ref{E:rec1}),
we have obtained that
\begin{equation}\label{inter:ynt}
Y_N(t)=U_N+M_{1,N}(t)+M_{2,N}(t)-V_{1,N}(t)
+\tilde V_{2,N}(t)-\frac{1}{2\hn^{d/2}}
\int_0^t\frac{d\langle Z_N\rangle _s}{Z_N^2(s)},
\end{equation}
where
$$
\tilde V_{2,N}(t)=
\left(\frac{\beta^2}{2\hat N^{3d/2}}\sum_{(i,j)\in
C_N}\qnd(i-j)+\frac{\beta^2\hat r\hat N^{d/2}}{2}\right)t
-\frac{\beta^2\gamma_0\hat N^{d\slash 2}t}{4}\left(1-r\right)^2.$$
Let us compute now the term $d\langle Z_N\rangle _s/Z_N^2(s)$:
according to Proposition \ref{P:exp}, $Z_N(t)$ is a continuous
semi-martingale, whose martingale part is
\begin{eqnarray*}
\hat M_N(t)&=&\frac{\beta}{\hat{N}^{d\slash
2}}\sum_{\sigma\in\Sigma_N}\sum_{(i,j)\in
                  C_N}\qn(i-j)\sigma_i\sigma_j\int_0^te^{-H_N(s)}dB_{i,j}(s)\\
                &&+\beta\hat r^{1\slash 2}\sum_{\sigma\in\Sigma_N}\sum_{i\in
                  C_N}\sigma_i\int_0^te^{-H_N(s)}dW_i(s).
\end{eqnarray*}
Hence,
\begin{eqnarray*}
\int_0^t\frac{1}{Z_N^2(s)}d\langle Z_N\rangle
_s&=&\int_0^t\frac{1}{Z_N^2(s)}d\langle\hat M_N\rangle _s\\
&=&\beta^2\hat r\sum_{i\in
  C_N}\int_0^t\rho_s(\sigma_i^1\sigma_i^2)ds
+\frac{\beta^2}{\hat N^d}\sum_{(i,j)\in
  C_N}\qnd(i-j)\int_0^t\rho_s\left(\sigma_i^1\sigma_j^1\sigma_i^2\sigma_j^2\right)ds.
\end{eqnarray*}
Recall now that $q^2(x)$ can be decomposed into
$q^2(x)=
\sum_{k\in\mathbb{Z}^d}\gamma_ke^{\imath \pi k\cdot  x}$. Thus,
\begin{multline*}
\int_0^t\frac{d\langle Z_N\rangle_s}{Z_N^2(s)}=\frac{\beta^2}{\hat
N^d}\sum_{(i,j)\in
  C_N}\sum_{k\in\mathbb{Z}^d}\gamma_k e^{\imath\pi k\cdot(i-j)\slash N}\int_0^t\rho_s\left(\sigma_i^1\sigma_j^1\sigma_i^2\sigma_j^2\right)ds\\
 +\beta^2\hat r\sum_{i\in
  C_N}\int_0^t\rho_s(\sigma_i^1\sigma_i^2)ds.
\end{multline*}
In order to simplify this last expression, we  will use the
elementary identity
\begin{align}\label{eq:modul-complex}
\left|\sum_{i\in C_N}z_i\right|^2=\sum_{i\in C_N}|z_i|^2+2\sum_{(i,j)\in
C_N}z_i\overline z_j,
\end{align}
valid for any family of complex numbers, and applied here to
$z_i=e^{\imath\pi k\cdot i/N}\sii^1\sii^2$. Recalling furthermore the 
definition (\ref{D: R_k}) of $R_k^{1,2}$, we end up with
\begin{multline*}
\int_0^t\frac{d\langle Z_N\rangle_s}{Z_N^2(s)} =
\frac{\hat N^d\beta^2}{2}
\sum_{k\in\mathbb{Z}_*^d}\gamma_k
\int_0^t\rho_s\left(\left|R_k^{1,2}\right|^2\right)ds
-\frac{\beta^2\gga t}{2}\\
+\frac{\beta^2\gamma_0\hat N^d}{2}\int_0^t\rho_s\left(\left(R^{1,2}\right)^2\right)ds
+\beta^2\hat N^{d}\hat r\int_0^t\rho_s(R^{1,2})ds.
\end{multline*}
Let us introduce now artificially the quantity $R^{1,2}_0$,
by writing $(R^{1,2})^2=(R^{1,2}_0)^2+2rR^{1,2}_0+r^2 $. This gives
\begin{multline}\label{E:rec2}
\int_0^t\frac{d\langle Z_N\rangle_s}{Z_N^2(s)} =
\frac{\hat N^d\beta^2}{2}\sum_{k\in\mathbb{Z}^d}\gamma_k\int_0^t\rho_s\left(\left|R_k^{1,2}\right|^2\right)ds-
\frac{\beta^2\gga t}{2}\\
-\frac{\beta^2\hat N^d\gamma_0r^2t}{2}
+2\beta^2\hat N^{d}\hat r\int_0^t\rho_s(R^{1,2})ds.
\end{multline}

\vspace{0.3cm}

Similarly to \cite{Ti}, let us notice that, since $X_i(s)\sim\cn(0,1-s)$,
a simple Gaussian integration by parts yields
\begin{equation}\label{gauss:ibp}
\E\lc \rho_s(\sigma_i)\frac{X_i(s)}{1-s}\rc=
\E\lc \partial_{X_i(s)}\rho_s(\sigma_i)\rc
=\beta
\hat r^{1\slash 2}\E\lc\rho_s(1-\sigma_i^1\sigma_i^2)\rc.
\end{equation}
This elementary consideration will induce us to add and subtract 
$\frac{\beta^2\hr}{\hat N^{d\slash 2}}\sum_i\rho_s(1-\sigma_i^1\sigma_i^2)$ to the expression
(\ref{inter:ynt}). By plugging moreover (\ref{E:rec2}), we get
\begin{equation}\label{inter:log-zn}
Y_N(t)=U_N+M_{1,N}(t)+M_{2,N}(t)-(V_{1,N}(t)-V_{2,N}(t))+\hat V_{3,N}(t),
\end{equation}
where
\begin{align*}\hat V_{3,N}(t)&=-\frac{\beta^2\hat r}{\hat
N^{d\slash 2}}\sum_{i\in
C_N}\int_0^t\rho_s(1-\sigma_i^1\sigma_i^2)ds
+\frac{\beta^2}{2\hat
N^{3d\slash 2}}\sum_{(i,j)\in
C_N}\qnd(i-j)t\\&+\frac{\beta^2\hat r}{2}\hat
N^{d\slash 2}t
-\frac{\hat
N^{d\slash 2}\beta^2}{4}\sum_{k\in\mathbb{Z}^d}\gamma_k
\int_0^t\rho_s\left(\left|R_k^{1,2}\right|^2\right)ds
+\frac{\beta^2 \gga t}{4\hat N^{d\slash 2}}
\\
&+\frac{\beta^2\hat N^{d\slash 2}\gamma_0r^2t}{4}
-\beta^2\hat N^{d\slash 2}\hat r\int_0^t\rho_s(R^{1,2})ds-\frac{\beta^2\gamma_0\hat N^{d\slash 2}t}{4}\left(1-r\right)^2.
\end{align*}
This last expression can be  simplified a little to give
\begin{multline*}
\hat V_{3,N}(t)=
\lc \frac{\gga}{\hn^{d/2}}
+\frac{2}{\hat N^{3d\slash 2}}\sum_{(i,j)\in
  C_N}\qnd(i-j)-\gamma_0\hat N^{d\slash 2}\rc
\frac{\beta^2 t}{4}\\
-
\frac{\hat
N^{d\slash 2}\beta^2}{4}
\sum_{k\in\mathbb{Z}^d}\gamma_k
\int_0^t\rho_s\left(\left|R_k^{1,2}\right|^2\right)ds,
\end{multline*}
that is
$\hat V_{3,N}(t)=V_{3,N}(t)$. By reporting this equality in
(\ref{inter:log-zn}), the proof is now complete.

\end{proof}

\vspace{0.3cm}

The last preliminary result we will need in order to establish 
our CLT is a self-averaging result for $R_k^{1,2}$ under 
the measure $\rho_t$ at a fixed value of the parameter $t\in\ou$.
\begin{proposition}\label{P:ref}
Under the hypotheses of Corollary \ref{cor1}, let $t\in\ou$. Then
\begin{align}
\sum_{k\in \mathbb{Z}^d}\gamma_k
\E\lc\rho_t(|R^{1,2}_k|^2)\rc=O\left(\frac{1}{N^{\hat d}}\right),
\end{align}
for $\hat d>d/2$, uniformly in $t$.
\end{proposition}

\begin{proof}
This is an easy elaboration of the computations leading to 
Proposition \ref{prop1}, by considerering the path
$v\in[0,t]\mapsto H_{N,t,v}(\si)$ defined by
\begin{equation*}
-H_{N,t,v}(\sigma)=-H_{N,t}^m(\sigma)+
\frac{\beta}{\hat{N}^{d\slash 2}}
\sum_{i\in C_N^m}B_{i,m}(v)\qn(i-m)\sigma_i\sigma_m
+\beta
\hat r^{1\slash 2} X_m(v)\sigma_m+h\sigma_m,
\end{equation*}
with obvious notations.

\end{proof}

\subsection{Proof of Theorem \ref{T:tcl}}

For an arbitrary $u\in\R$, we will try to control
$\Phi_{N,u}(t)\equiv\E [e^{\imath uY_N(t)}]$. To this purpose, we will
apply It\^o's formula to the complex valued $C_b^2$ function 
$x\mapsto e^{\imath ux}$. We obtain, for any $t\in\ou$,
\begin{align*}
e^{\imath u Y_N(t)}= D_{1,N}+\sum_{m=2}^8D_{m,N}(t),
\end{align*}
where
\begin{align*}
D_{1,N}&=e^{\imath u U_N}\\
D_{2,N}(t)&=\frac{\imath u\beta}{\hat N^{d}}\sum_{(i,j)\in
C_N}\qn(i-j)\int_0^te^{\imath u Y_N(s)}\rho_s(\sigma_i\sigma_j)dB_{i,j}\\
D_{3,N}(t)&=\frac{\imath u\beta\hat r^{1\slash 2}}{\hat N^{d\slash
2}}\sum_{i\in C_N}\int_0^te^{\imath u Y_N(s)}\rho_s(\sigma_i)dW_i\\
 D_{4,N}(t)&=-\frac{\imath
u\beta\hat r^{1\slash 2}}{\hat N^{d\slash 2}}\sum_{i\in
C_N}\int_0^te^{\imath u Y_N(s)}\rho_s(\sigma_i)\frac{X_i(s)}{1-s}ds   \\
D_{5,N}(t)&=\frac{\imath u\beta^2\hat r}{\hat N^{d\slash
2}}\sum_{i\in C_N}\int_0^te^{\imath u
Y_N(s)}\rho_s(1-\sigma_i^1\sigma_i^2)ds,   
\end{align*}
and
\begin{multline*}
 D_{6,N}(t)=
\frac{\imath u\beta^2}{4}\lc  \frac{\gga}{\hn^{d/2}}
+\frac{2}{\hat N^{3d\slash 2}}\sum_{(i,j)\in C_N}\qnd(i-j)
- \gamma_0\hn^{d/2}\rc
\int^t_0e^{\imath uY_N(s)}ds\\
-
\frac{\imath u\hat N^{d\slash 2}\beta^2}{4}
\sum_{k\in\mathbb{Z}^d}\gamma_k \int_0^te^{\imath uY_N(s)}
\rho_s\left(\left|R_k^{1,2}\right|^2\right)ds.
\end{multline*}
The terms $D_{7,N}(t)$ and $D_{8,N}(t)$are the ones associated with the 
quadratic variation process of $Y_N(t)$, that is:
\begin{align*}
 D_{7,N}(t)&=-\frac{u^2\beta^2}{2\hat
N^{2d}}\sum_{(i,j)\in C_N}\qnd(i-j)\int_0^te^{\imath u
Y_N(s)}\rho_s(\sigma_i^1\sigma_j^1\sigma^2_i\sigma_j^2)ds\\
&=-\frac{u^2\beta^2}{2\hat N^{2d}}\sum_{(i,j)\in C_N}\sum_{k\in
\mathbb{Z}^d}\gamma_ke^{\imath \pi k\cdot(i-j)/N}\int_0^te^{\imath
uY_N(s)}\rho_s(\sigma_i^1\sigma_j^1\sigma^2_i\sigma_j^2)ds\\
 D_{8,N}(t)&=-\frac{u^2\beta^2\hat r}{2\hat N^d}\sum_{i\in
C_N}\int_0^te^{\imath u
Y_N(s)}\rho_s(\sigma_i^1\sigma_i^2)ds=-\frac{u^2\beta^2\hat
r}{2}\int_0^te^{\imath u Y_N(s)}\rho_s(R^{1,2})ds. 
\end{align*}

\vspace{0.3cm}

Let us find some estimates for the expected value of all the terms 
we have obtained in this decomposition. First of all, the usual
central limit theorem for {\small IID} random variables yields
\begin{align}
\E [D_{1,N}(t)]=e^{-\hat\tau^2u^2\slash 2}
+O\left(\frac{1}{N^{d\slash 2}}\right),
\end{align}
where $\hat\tau$ is defined at (\ref{def:tau}).
Furthermore, the terms $D_{2,N}(t)$ and $D_{3,N}(t)$ are obviously
of zero mean.

\vspace{0.3cm}

In order to control $D_{4,N}(t)$, let us perform again the Gausssian
integration by parts (\ref{gauss:ibp}), which can be read here as
\begin{eqnarray*}
\E [D_{4,N}(t)]&=&-\frac{\imath u\beta \hat r^{1\slash 2}}{\hat
N^{d \slash 2}}\sum_{i\in C_N}\int_0^t\E 
\lc \partial_{X_i(s)}(e^{\imath u Y_N(s)}\rho_s(\sigma_i))\rc ds\\
&=&\frac{u^2\beta\hat r^{1\slash 2}}{\hat
N^{d\slash 2}}\sum_{i\in
C_N}\int_0^t\E\lc \partial_{X_i(s)}(Y_N(s)) e^{\imath u
Y_N(s)}\rho_s(\sigma_i)\rc ds-\E [D_{5,N}(t)].
\end{eqnarray*}
Furthermore,
\begin{align}
\partial_{X_i(s)}(Y_N(s))=\frac{\partial_{X_i(s)}(Z_N(s))}{N^{d\slash
2}Z_N(s)}=\frac{\beta\hat r^{1\slash 2}}{N^{d\slash
2}}\rho_s(\sigma_i),\nonumber
\end{align}
and therefore, we obtain
\begin{align}
\E [D_{4,N}(t)+D_{5,N}(t)]&=u^2\beta^2\hat
r\int^t_0\E \lc e^{\imath uY_N(s)}\rho_s(R^{1,2})\rc ds\nonumber\\
&=u^2\beta^2\gamma_0
r^2\int_0^t\Phi_{N,u}(s)ds+O\left(\frac{1}{N^{\hat d}}\right),
\end{align}
owing to Proposition \ref{P:ref}.

\vspace{0.3cm}

Let us turn now to the estimation of $D_{6,N}(t)$: notice that, under
Hypothesis \ref{hyp:q} and \ref{hyp:q2}, it is easily checked that
$$
\left|\frac{1}{\hat N^{2d}}\sum_{(i,j)\in
    C_N}\qnd(i-j)-\frac{\gamma_0}{2}\right|=O\left(\frac{1}{\hat
    N^{\hat d}}\right),
$$
which together with a direct application of Proposition \ref{P:ref},
shows that
\begin{align*}
\E [D_{6,N}(t)]=O\left(\frac{1}{\hat N^\epsilon}\right),
\end{align*}
where $\epsilon=\hat d -d\slash 2$.

\vspace{0.3cm}

As far as  $D_{7,N}(t)$ is concerned, notice 
that by relation (\ref{eq:modul-complex}), we have
\begin{align}
 D_{7,N}(t)&=-\frac{u^2\beta^2}{4}\sum_{k\in
\mathbb{Z}_*^d}\gamma_k\int_0^te^{\imath u
Y_N(s)}\rho_s\left(|R^{1,2}_k|^2\right)ds\nonumber\\
&-\frac{u^2\beta^2\gamma_0}{4}\int_0^te^{\imath u
Y_N(s)}\rho_s\left(|R^{1,2}|^2\right)ds+\frac{u^2\beta^2\gga}{4\hat N^d}
\int_0^te^{\imath u Y_N(s)}ds,
\end{align}
and thanks  to proposition \ref{P:ref}, we get
\begin{align}
 D_{7,N}(t)&=-\frac{u^2\beta^2\gamma_0r^2}{4}\int_0^te^{\imath u
Y_N(s)}ds+ O\left(\frac{1}{N^{\hat d}}\right),
\end{align}
and thus
\begin{align*}
\E [D_{7,N}(t)]&=-\frac{u^2\beta^2\gamma_0
r^2}{4}\int_0^t\Phi_{N,u}(s)ds+O\left(\frac{1}{N^{\hat d}}\right).
\end{align*}
Eventually, in a similar way we have
\begin{align*}
\E \lc D_{8,N}(t)\rc=-\frac{u^2\beta^2\gamma_0
r^2}{2}\int_0^t\Phi_{N,u}(s)ds + O\left(\frac{1}{\hat N^{\hat d}}\right).
\end{align*}

\vspace{0.3cm}

Putting together the previous estimates on 
$\E[D_{1,N}(t)],\ldots,\E[D_{8,N}(t)]$, we have finally:
$$
\psi_{N,u}(t)=e^{-\frac{\nu^2u^2}{2}}
+ \frac{(u\beta \sqrt{\gamma_0}r)^{2}}{4}\int_0^t  \psi_{N,u}(s) ds+
\hat R_{N,u}(t),
$$
with
$$
\left|\hat R_{N,u}(t)  \right|\le \frac{\kappa}{N^{\epsilon}},
$$
which ends the proof by a Gronwall type argument.

\hfill
$\Box$

\vspace{0.5cm}

\noindent
{\bf Acknowledgement:}
We would like to thank Frederi Viens for letting us take 
advantage of some old computations for the simple limit 
of the free energy of the localized model.

\end{document}